# FREUD'S EQUATIONS FOR ORTHOGONAL POLYNOMIALS AS DISCRETE PAINLEVÉ EQUATIONS.

Alphonse P. Magnus

ABSTRACT. We consider orthogonal polynomials $p_n$ with respect to an exponential weight function $w(x) = \exp(-P(x))$. The related equations for the recurrence coefficients have been explored by many people, starting essentially with Laguerre [49], in order to study special continued fractions, recurrence relations, and various asymptotic expansions (G. Freud's contribution [28,56]).

Most striking example is $n = 2tw_n + w_n(w_{n+1} + w_n + w_{n-1})$ for the recurrence coefficients $p_{n+1} = xp_n - w_n p_{n-1}$ of the orthogonal polynomials related to the weight $w(x) = \exp(-4(tx^2 + x^4))$ (notations of [26, pp.34–36]). This example appears in practically all the references below. The connection with discrete Painlevé equations is described here.

## 1. CONSTRUCTION OF ORTHOGONAL POLYNOMIALS RECURRENCE COEFFICIENTS.

Consider the set $\{p_n\}_0^\infty$ of orthonormal polynomials with respect to a weight function $w$ on (a part of) the real line:

$$\int_{-\infty}^{\infty} p_n(x) p_m(x) w(x)\, dx = \delta_{n,m}, \qquad n, m = 0, 1, \ldots \qquad (1)$$

The very useful *recurrence formula* has the form $a_1 p_1(x) = (x - b_0) p_0$,

$$a_{n+1} p_{n+1}(x) = (x - b_n) p_n(x) - a_n p_{n-1}(x), \quad n = 1, 2, \ldots \qquad (2)$$

The connection between the nonnegative integrable function $w$ and the real sequences $\{a_n\}_1^\infty, \{b_n\}_0^\infty$ is of the widest interest. It is investigated in combinatorics [9,10,52], asymptotic analysis [11,12,28,30,53–55,61,77–79], numerical analysis [15, 21, 31, 32, 51, 57], and, of course, Lax-Painlevé-Toda theory (all the other references, excepting *perhaps* Chaucer).

However, the connection may seem quite elementary and explicit:

1. From $w$ to $a_n$ and $b_n$: let $\mu_k = \int_{-\infty}^{\infty} x^k w(x)\, dx$, $k = 0, 1, \ldots$ be the moments of $w$, then

$$a_n^2 = \frac{H_{n-1} H_{n+1}}{(H_n)^2}, n = 1, 2, \ldots \qquad b_n = \frac{\widetilde{H}_{n+1}}{H_{n+1}} - \frac{\widetilde{H}_n}{H_n}, n = 0, 1, \ldots \qquad (3)$$

where $H_n$ and $\widetilde{H}_n$ are the determinants

$$H_n = \begin{vmatrix} \mu_0 & \cdots & \mu_{n-1} \\ \vdots & & \vdots \\ \mu_{n-1} & \cdots & \mu_{2n-2} \end{vmatrix} ; \widetilde{H}_n = \begin{vmatrix} \mu_0 & \cdots & \mu_{n-2} & \mu_n \\ \vdots & & \vdots & \vdots \\ \mu_{n-1} & \cdots & \mu_{2n-3} & \mu_{2n-1} \end{vmatrix}, \qquad (4)$$

$n = 0, 1, \ldots, H_0 = 1, \widetilde{H}_0 = 0$.

2. From the $a_n$'s and $b_n$'s to $w$: consider the Jacobi matrix

$$\boldsymbol{J} = \begin{bmatrix} b_0 & a_1 & & & \\ a_1 & b_1 & a_2 & & \\ & a_2 & b_2 & a_3 & \\ & & \ddots & \ddots & \ddots \end{bmatrix}$$

then, remarking that (2) may be written $xp_n(x) = a_n p_{n-1}(x) + b_n p_n(x) + a_{n+1} p_{n+1}(x)$, (Remark that the coefficients appear now with positive signs, a very elementary fact, but always able to puzzle weak souls), or $x\boldsymbol{p}(x) = \boldsymbol{J}\boldsymbol{p}(x)$, where $\boldsymbol{p}(x)$ is the column (infinite) vector of $p_0(x), p_1(x), \ldots$, so that the expansion of $x^k p_n(x)$ in the basis $\{p_0, p_1, \ldots\}$ is $x^k p_n(x) = \sum_\ell (\boldsymbol{J}^k)_{n,\ell} p_\ell(x)$. As we are dealing with a basis of orthonormal polynomials, this means that $(\boldsymbol{J}^k)_{n,\ell}$ is the $\ell^{\text{th}}$ Fourier coefficient of $x^k p_n(x)$, (where row and column indexes start at zero), i.e., $(\boldsymbol{J}^k)_{n,\ell} = \int_{-\infty}^{\infty} x^k p_n(x) p_\ell(x) w(x)\, dx$. For any polynomial $p$, one has

$$(p(\boldsymbol{J}))_{n,\ell} = \int_{-\infty}^{\infty} p(x) p_n(x) p_\ell(x) w(x)\, dx. \qquad (5)$$

In particular, $(\boldsymbol{J}^k)_{0,0} = \mu_k/\mu_0$.

Direct numerical use of (3) is almost always unsatisfactory, for stability reasons. Together with various ways to cope with the numerical stability problem [31,32], compact formulas or equations for the recurrence coefficients $a_n$ and $b_n$ in special cases have been sought. Of course, quite a number of exact solutions are known (appendix of [14]) and new ones are steadily discovered (Askey-Wilson polynomials and other instances of the new "$q$−calculus" [4, 35, 36, 47]). The trend explored here is essentially related to weight functions satisfying simple differential equations. Main innovators were Laguerre in 1885 [49], Shohat in 1939 [77], and Freud in 1976 [28]. Their contributions are now described in reverse order. On Géza Freud (1922-1979), see [29,68,69].

## 2. Definition of Freud's equations.

In order to reduce the technical contents of what will follow, only *even* weights: $\forall$ real $x$, $w(-x) = w(x)$ will be used. Then, the orthogonal polynomials are even or odd according to their degrees; equivalently, all the $b_n$'s in (2), in (3), and in the Jacobi matrix $\boldsymbol{J}$ vanish. Let $\boldsymbol{a}$ be the sequence of coefficients $\{a_1, a_2, \ldots\}$.

**Theorem 1.** *Let $P$ be the real even polynomial $P(x) = c_0 x^{2m} + c_1 x^{2m-2} + \cdots + c_m$, with $c_0 > 0$, and let $a_1, a_2, \ldots$ be the recurrence coefficients (2) of* *the orthonormal polynomials with respect to the weight* $w(x) = \exp(-P(x))$ *on the whole real line. Then, the* **Freud's equations**

$$F_n(\boldsymbol{a}) := a_n \left(P'(\boldsymbol{J})\right)_{n,n-1} = n, \quad n = 1, 2, \ldots \tag{6}$$

hold, where $(\ )_{n,n-1}$ means the element at row $n$ and column $n-1$ of the (infinite) matrix $P'(\boldsymbol{J})$. A finite number of computations is involved in each of these elements. The row and column indexes start at 0. We could as well take the indexes $(n-1, n)$, as $\boldsymbol{J}$ and any polynomial function of $\boldsymbol{J}$ are symmetric matrices.

**Remark.** One has

$$F_1(\boldsymbol{a}) = (\boldsymbol{J} P'(\boldsymbol{J}))_{0,0}. \tag{7}$$

Indeed, the first element of the first row of $\boldsymbol{J}P'(\boldsymbol{J})$ is $a_1$ times the element $(1, 0)$ of $P'(\boldsymbol{J})$, i.e., $F_1(\boldsymbol{a})$.

**Examples.**

$$w(x) = \exp(-x^2) \quad F_n(\boldsymbol{a}) = 2a_n^2, \tag{8}$$

$$w(x) = \exp(-\alpha x^4 - \beta x^2) \quad F_n(\boldsymbol{a}) = 4\alpha a_n^2 (a_{n-1}^2 + a_n^2 + a_{n+1}^2) + 2\beta a_n^2, \tag{9}$$

$$w(x) = \exp(-x^6) \quad F_n(\boldsymbol{a}) = 6a_n^2(a_{n-2}^2 a_{n-1}^2 + a_{n-1}^4 + 2a_{n-1}^2 a_n^2 + a_n^4$$
$$+ a_{n-1}^2 a_{n+1}^2 + 2a_n^2 a_{n+1}^2 + a_{n+1}^4 + a_{n+1}^2 a_{n+2}^2), \tag{10}$$

The elementary Hermite polynomials case is of course immediately recovered in (8). (9) and (10) were used by Freud [28] in investigations on asymptotic behaviour, see [11, 53, 56, 57, 67, 69] for more. The rich connection of (9) with discrete and continuous Painlevé theory will be recalled later on. For a much more general setting, see [1, Theor. 4.1].

**Proof of** (6). We consider two different ways to write the integral $I_n = \int_{-\infty}^{\infty} p'_n(x) p_{n-1}(x) \exp(-P(x))\, dx$. First, let $p_{n-1}(x) = \pi_{n-1} x^{n-1} + \ldots$ and $p_n(x) = \pi_n x^n + \ldots$. From (2) (with $n-1$ instead of $n$), $\pi_n = \pi_{n-1}/a_n$; so, $p'_n(x) = n\pi_n x^{n-1} + \ldots = \dfrac{n}{a_n} p_{n-1}(x) +$ a polynomial of degree $\leqslant n-2$. As $p_{n-1}$ is orthogonal to all polynomials of degree $\leqslant n-2$, what remains is $I_n = \dfrac{n}{a_n}$, as the orthonormal polynomials have a unit square integral. Next, we perform an integration by parts on $I_n = -\int_{-\infty}^{\infty} p_n(x) \left[p_{n-1}(x) \exp(-P(x))\right]'\, dx = \int_{-\infty}^{\infty} P'(x) p_n(x) p_{n-1}(x) \exp(-P(x))\, dx$, using the orthogonality of $p_n$ and $p'_{n-1}$ of degree $< n-1 < n$, and the latter integral is $(P'(\boldsymbol{J}))_{n,n-1}$, according to the spectral representation (5). □

### 3. Freud's equations software.

The interested reader (if there is still one left) will find in the "software" part of http://www.math.ucl.ac.be/~magnus/ three FORTRAN programs related to Freud's equations. `freud1.f` asks for an exponent $m$ and puts data associated to the Freud's equations for $x^2, x^4, \ldots, x^{2m}$ in the file `freud00m.dat`. `freud2.f` produces a readable (sort of) listing out of such a data file. Finally, `freud3.f` reads the coefficients of a polynomial $P$ and computes in a stable way a sequence of positive recurrence coefficients $a_1, a_2, \ldots$ associated to the weight $\exp(-P(x))$ on the whole real line. Actually, weight functions $|x|^\alpha \exp(-P(x))$, with $P$ even, on the whole real line, and $x^\beta \exp(-Q(x))$ on the *positive* real line are considered simultaneously ($P(x) = Q(x^2)$).

### 4. Solution of Freud's equations.

We started from the *solution* (3), and built the *equation* (6) afterwards, when the $a_n$'s are recurrence coefficients of orthogonal polynomials associated to the weight $\exp(-P(x))$.

We try now to discuss the general solution of (6). With $P$ even of degree $2m$, the $n^{\text{th}}$ equation of (6) involves $a_{n-m+1}, \ldots, a_{n+m-1}$ (see [56, 57] or the examples (8) (9), (10)), so that the general solution depends on $2m - 2$ arbitrary constants, for instance $a_{n_0-m+1}, \ldots, a_{n_0+m-2}$. I shall only solve the case $a_{-m+1} = \cdots = a_0 = 0$, and investigate how the solution depends on the $m - 1$ initial data $a_1, \ldots, a_{m-1}$:

**Theorem 2.** *With $P(x) = c_0 x^{2m} + c_1 x^{2m-2} + \cdots + c_m$, the solution of (6) with $a_{-m+1} = \cdots = a_0 = 0$ is (3), with (4), where the $\mu_j$'s of (4) satisfy the linear recurrence relation of order $2m$*

$$2mc_0 \mu_{2n+2m} + (2m-2)c_1 \mu_{2n+2m-2} + \cdots + 2c_{m-1} \mu_{2n+2} = (2n+1)\mu_{2n}, \quad (11)$$

*for $n = 0, 1, \ldots$, and $\mu_{2j+1} = 0$. $\mu_2, \mu_4, \ldots, \mu_{2m-2}$ are given by $\mu_{2j}/\mu_0 = (\boldsymbol{J}^{2j})_{0,0}$.*

The following lemma will be used:

**Lemma 1.** *Let $\boldsymbol{a}^{(k)} = \{a_{n,k}\}_{n=1}^\infty$, $k = 0, 1, \ldots$, be the solutions of the **quotient-difference** equations*

$$\begin{aligned} a_{n-1,k+1}^2 + a_{n,k+1}^2 &= a_{n,k}^2 + a_{n+1,k}^2, & \text{if } n \text{ is odd}, \\ a_{n-1,k+1} a_{n,k+1} &= a_{n,k} a_{n+1,k}, & \text{if } n \text{ is even}, \end{aligned} \quad (12)$$

*then, $F_n(\boldsymbol{a}^{(k)})$ satisfies*

$$\frac{a_{n-1,k+1}}{a_{n,k}} \left( F_n(\boldsymbol{a}^{(k+1)}) - F_n(\boldsymbol{a}^{(k)}) \right) = \frac{a_{n,k}}{a_{n-1,k+1}} \left( F_{n+1}(\boldsymbol{a}^{(k)}) - F_{n-1}(\boldsymbol{a}^{(k+1)}) \right),$$

$$n \text{ even},$$

$$F_n(\boldsymbol{a}^{(k+1)}) - F_n(\boldsymbol{a}^{(k)}) = F_{n+1}(\boldsymbol{a}^{(k)}) - F_{n-1}(\boldsymbol{a}^{(k+1)}), \quad n \text{ odd}. \quad (13)$$

where $F_n(\boldsymbol{a}^{(k)}) = a_{n,k}\,(P'(\boldsymbol{J}_k))_{n,n-1}$, and where $\boldsymbol{J}_k$ is the Jacobi matrix

$$\boldsymbol{J}_k = \begin{bmatrix} 0 & a_{1,k} & & & \\ a_{1,k} & 0 & a_{2,k} & & \\ & a_{2,k} & 0 & a_{3,k} & \\ & & \ddots & \ddots & \ddots \end{bmatrix}$$

The quotient-difference equations are well known in orthogonal polynomials identities investigations, see [74, 75] for recent examples. The notation chosen here is motivated by the interpretation of $a_{n,k}$ as recurrence coefficient of orthogonal polynomials related to the weight $x^{2k}w(x)$, but we do not need this interpretation in order to establish the theorem and the lemma. Actually, the theorem 2 is meant to recover an interpretation of the solutions of (6) in terms of orthogonal polynomials, but arbitrary initial data $a_1, a_2, \ldots, a_{m-1}$ will yield formal orthogonal polynomials through formal moments $\mu_j$ which, from (11), can still be written $\mu_j = \int_S x^j \exp(-P(x))\,dx$, but where the support $S$ is now a system of arcs in the complex plane, see [66, "Laplace's method"], [62, 63].

Remark that, with $a_{n,0} = a_n$ and $a_{0,k} = 0$, all the $a_{n,k}$'s are completely determined by (12) as functions of $a_1, a_2, \ldots$

**Proof of Lemma 1.** The quotient-difference equations (12) can be written $\boldsymbol{J}_k \boldsymbol{L}_k = \boldsymbol{L}_k \boldsymbol{J}_{k+1}$, where $\boldsymbol{L}_k$ is the infinite lower triangular matrix

$$\boldsymbol{L}_k = \begin{bmatrix} a_{1,k} & & & & & \\ 0 & a_{1,k+1} & & & & \\ a_{2,k} & 0 & a_{3,k} & & & \\ & a_{2,k+1} & 0 & a_{3,k+1} & & \\ & & a_{4,k} & 0 & a_{5,k} & \\ & & & \ddots & \ddots & \ddots \end{bmatrix}$$

One finds also $\boldsymbol{J}_k^2 = \boldsymbol{L}_k \boldsymbol{R}_k$, $\boldsymbol{J}_{k+1}^2 = \boldsymbol{R}_k \boldsymbol{L}_k$, where $\boldsymbol{R}_k$ is the transposed of $\boldsymbol{L}_k$ (the famous $LR$ relations). So, $\boldsymbol{J}_k^3 \boldsymbol{L}_k = \boldsymbol{L}_k \boldsymbol{R}_k \boldsymbol{L}_k \boldsymbol{J}_{k+1} = \boldsymbol{L}_k \boldsymbol{J}_{k+1}^3$, ..., $\boldsymbol{J}_k^{2p+1} \boldsymbol{L}_k = \boldsymbol{L}_k \boldsymbol{R}_k \boldsymbol{L}_k \boldsymbol{J}_{k+1}^{2p-1} = \boldsymbol{L}_k \boldsymbol{J}_{k+1}^{2p+1}$, for any odd power, so $P'(\boldsymbol{J}_k)\boldsymbol{L}_k = \boldsymbol{L}_k P'(\boldsymbol{J}_{k+1})$, and we look at the element of indexes $(n, n-1)$:

$$(P'(\boldsymbol{J}_k))_{n,n-1}(\boldsymbol{L}_k)_{n-1,n-1} + (P'(\boldsymbol{J}_k))_{n,n+1}(\boldsymbol{L}_k)_{n+1,n-1} = \\ (\boldsymbol{L}_k)_{n,n-2}(P'(\boldsymbol{J}_{k+1}))_{n-2,n-1} + (\boldsymbol{L}_k)_{n,n}(P'(\boldsymbol{J}_{k+1}))_{n,n-1},$$

or $F_n(\boldsymbol{a}^k)\dfrac{(\boldsymbol{L}_k)_{n-1,n-1}}{a_{n,k}} + F_{n+1}(\boldsymbol{a}^k)\dfrac{(\boldsymbol{L}_k)_{n+1,n-1}}{a_{n+1,k}} =$

$$\dfrac{(\boldsymbol{L}_k)_{n,n-2}}{a_{n-1,k+1}}F_{n-1}(\boldsymbol{a}^{k+1}) + \dfrac{(\boldsymbol{L}_k)_{n,n}}{a_{n,k+1}}F_n(\boldsymbol{a}^{k+1}),$$

leading to (13). □

**Proof of Theorem 2.** From (13), with (6), $(\boldsymbol{a}^{(0)} = \boldsymbol{a})$, and $F_0(\boldsymbol{a}^{(k)}) = 0$ for all $k$, one finds by induction $F_n(\boldsymbol{a}^{(k)}) = n + k(1 - (-1)^n)$, so, using (7), $F_1(\boldsymbol{a}^{(k)}) = (\boldsymbol{J}_k P'(\boldsymbol{J}_k))_{0,0} = 2k+1$, $k = 0, 1, \ldots$ Let $\mu_{2k}$ = constant $\times a_1^2 a_{1,1}^2 a_{1,2}^2 \cdots a_{1,k-1}^2$ and $\mu_{2k+1} = 0$. From the $LR$ relations, $(\boldsymbol{J}_k^{2p})_{0,0} = (\boldsymbol{L}_k \boldsymbol{L}_{k+1} \cdots \boldsymbol{L}_{k+p} \boldsymbol{R}_{k+p} \cdots \boldsymbol{R}_{k+1} \boldsymbol{R}_k)_{0,0} = a_{1,k}^2 a_{1,k+1}^2 \cdots a_{1,k+p-1}^2 = \mu_{2k+2p}/\mu_{2k}$.

Finally, $2k + 1 = (\boldsymbol{J}_k P'(\boldsymbol{J}_k))_{0,0} = 2mc_0(\boldsymbol{J}_k^{2m})_{0,0} + (2m-2)c_1(\boldsymbol{J}_k^{2m-2})_{0,0} + \cdots + 2c_{m-1}(\boldsymbol{J}_k^2)_{0,0} = 2mc_0\mu_{2k+2m}/\mu_{2k} + (2m-2)c_1\mu_{2k+2m-2}/\mu_{2k} + \cdots + 2c_{m-1}\mu_{2k+2}/\mu_{2k}$, whence (11), and the determinant ratios (3) follow from (12): $a_{n,k}^2 = H_{n-1,k}H_{n+1,k}/(H_{n,k})^2$ where $H_{n,k}$ is the Hankel determinant of rows $\mu_{2k}, \mu_{2k+1}(=0), \ldots, \mu_{2k+n-1}$; $\mu_{2k+1}, \ldots, \mu_{2k+n}$, etc., up to $\mu_{2k+n-1}, \ldots, \mu_{2k+2n-2}$.    □

The *isomonodromy approach* looks at the linear differential equation satisfied by an orthogonal polynomial $p_n$ [5, 16, 17, 37, 38] and produces isomonodromy identities of interest here. For instance, working (9) leads to $a_n^2 = (4\alpha)^{-1/2}\mathrm{P}_{\mathrm{IV}}((4\alpha)^{-1/2}\beta; -n/2, -n^2/4)$ (Kitaev, [23–26,43]), where $\mathrm{P}_{\mathrm{IV}}(t; A, B)$ is the solution of the fourth Painlevé equation

$$\ddot{y} = \frac{\dot{y}^2}{2y} + \frac{3y^3}{2} + 4ty^2 + 2(t^2 - A)y + \frac{B}{y}$$

which remains $O(t^{-1})$ when $t \to +\infty$ [60, pp. 231-232] (the importance and relevance of this latter condition is not yet quite clear). See [1–3, 16–20, 22–26, 39, 43, 44, 46] for more on such connections.

## 5. Freud's equations as discrete Painlevé equations.

As far as I understand the matter, discrete Painlevé equations were first designed as clever discretisations of the genuine Painlevé equations, so to keep interesting features of these equations, especially integrability. For instance, (9) when $\beta \neq 0$ is a discretisation of the Painlevé-I equation [23–26, 33, 34, 43]. Then, these features are examined on various discrete equations not necessarily associated to Painlevé equations. So, for instance, [33, p.350] encounter (9) with $\beta = 0$ as a discrete equation which is no more linked to the discretisation of a continuous Painlevé equation (and the authors of [33] call (9) a discrete Painlevé-0 equation when $\beta = 0$).

If we want to check that the Freud's equations are a valuable instance of discrete Painlevé equation, the following features are of interest, according to experts [48]:

**5.1. Analyticity.** Each $a_{n+m-1}^2$ is a rational function of preceding elements, so a meromorphic function of initial conditions.

**5.2. Reversibility.** There is no way to distinguish past and future, (6) is left unchanged when $a_{n+1}, a_{n+2}, \ldots$ are permuted with $a_{n-1}, a_{n-2}, \ldots$ Indeed ([56, p.369]), (6) is a sum of terms $a_n^2 a_{n+i_1}^2 a_{n+i_2}^2 \cdots a_{n+i_p}^2$, provided $0 \leqslant i_1 + 1 \leqslant i_2 + 2 \leqslant \cdots \leqslant i_p + p \leqslant p + 1$. Then, $\{-i_1, -i_2, \ldots, -i_p\}$ satisfies the same conditions and is therefore present too: $0 \leqslant -i_p + 1 \leqslant -i_{p-1} + 2 \leqslant \cdots \leqslant -i_1 + p \leqslant p + 1$.

**5.3. Symmetry.** An even stronger property is the following:

**Theorem 3.** *The matrix of derivatives* $\left(\dfrac{\partial F_n(\boldsymbol{a})}{\partial \log a_m}\right)_{n,m=1,2,\ldots}$ *is symmetric.*

Proof: see [56, 57].

**5.4. Integrability.** The existence of the formula (3) according to Theorem 2 shows how (6) can be solved. The similarity of (3) with known solutions of discrete Painlevé equations is also striking [42, 45]. The structure of the solution leads also to

**5.5. Singularity confinement.** The movable poles property of continuous equations is replaced by the following [34]: if initial values are such that some $a_n$, say $a_{n_0}$ becomes very large, only a finite number of neighbours $a_{n_0+1}, \ldots, a_{n_0+p}$ should be liable to be very large too; moreover, $a_{n_0+p+1}$, etc. should be continuous functions of $a_{n_0-1}, a_{n_0-2}$ etc. The formula (3) shows that the singularities of $a_n$ are the zeros of $H_n$ which are determinants built with solutions of (11), therefore continuous functions of the initial data.

## 6. Generalizations of Freud's equations: semiclassical orthogonal polynomials, etc.

Let the weight function satisfy the differential equation $w'/w = 2V/W$, where $V$ and $W$ are polynomials. The related orthogonal polynomials on a support $S$ are then called *semi-classical* [6–8, 40, 41, 64, 65, 72] (classical: degrees of $V$ and $W \leqslant 1$ and 2). One can still build equations for the recurrence coefficients by working out the integrals $I_n = \int_S W(x) p'_n(x) p_{n-1}(x) w(x) \, dx$, and $J_n = \int_S W(x) p'_n(x) \, p_n(x) w(x) \, dx$. Many equivalent forms can be found, and it not yet clear to know which one is the most convenient. If we manage to have $W(x)w(x) = 0$ at the endpoints of the support $S$, explicit formulas follow [77] (see also formula (4.5) of [1]). For instance, with $w(x) = \exp(-Ax^2)$ on $S = [-a, a]$ (see [15] for other equations), one takes $W(x) = x^2 - a^2$ and $V(x) = -Ax(x^2 - a^2)$, uses $p'_n = np_{n-1}/a_n + [2(a_1^2 + \cdots + a_{n-1}^2) - na_{n-1}^2]p_{n-3}/(a_{n-2}a_{n-1}a_n) + \cdots$ from (2), to find $2(a_1^2 + \cdots + a_{n-1}^2) + (2n+1)a_n^2 - na^2 - 2Aa_n^2(a_{n-1}^2 + a_n^2 + a_{n+1}^2 - a) = 0$, or $(2n+1)x_{n+1} - (2n-1)x_n - na^2 - 2A(x_{n+1} - x_n)(x_{n+2} - x_{n-1} - a) = 0$, where $x_n = a_1^2 + \cdots + a_{n-1}^2$.

It is not clear how to recover reversibility and symmetry. For an even much more nasty case, see [61].

But the root of the matter lies even deeper, as already recognized by Laguerre in 1885 [49]! The point is to derive equations for the coefficients in the Jacobi continued fraction $f(z) = \cfrac{1}{z - b_0 - \cfrac{a_1^2}{z - b_1 - \cdots}}$ from a differential equation $Wf' = 2Vf + U$ with polynomial coefficients (same $W$ and $V$ as above: orthogonal polynomials lead to this continued fraction through $f(z) = \int_S (z-x)^{-1} w(x) \, dx$). See [30, 72] for this technique.

The continued fraction approach leads readily to the *Riccati* extension $Wf' = Rf^2 + 2Vf + U$ [20,55], as each continued fraction

$$f_k(z) = \cfrac{1}{z - b_k - \cfrac{a_{k+1}^2}{z - b_{k+1} - \cdots}}$$

is found to satisfy a Riccati equation $Wf'_k = -a_k^2 \Theta_{k-1} f_k^2 + 2\Omega_k f_k - \Theta_k$ too, where $\Theta_k$ and $\Omega_k$ are polynomials of fixed degrees (i.e., independent. of $k$). Equations for the recurrence coefficients $a_k$'s and $b_k$'s follow from recurrence relations for the $\Theta_k$'s and $\Omega_k$'s (use $f_0(z) = f(z)$ and $f_k(z) = 1/[z - b_k - a_{k+1}^2 f_{k+1}(z)]$).

Finally, the differential operator can be replaced by a suitable difference operator, amounting to $W(z)\dfrac{f(\varphi_2(z)) - f(\varphi_1(z))}{\varphi_2(z) - \varphi_1(z)} = R(z)f(\varphi_1(z))f(\varphi_2(z)) + V(z)[f(\varphi_1(z)) + f(\varphi_2(z))] + U(z)$, where $\varphi_1(z)$ and $\varphi_2(z)$ are the two roots of a quadratic equation $A\varphi^2(z) + 2Bz\varphi(z) + Cz^2 + 2D\varphi(z) + 2Ez + F = 0$, see [36,58].

### Acknowledgements.

Thanks to the Conference organizers, and to Darwin's College hospitality, according to the highest quality standards since more than 600 years:

| | |
|---|---|
| 26: In felaweshipe, and pilgrimes were they alle, | In fellowship, and pilgrims were they all |
| 27: That toward caunterbury wolden ryde. | That toward Canterbury town would ride. |
| 28: The chambres and the stables weren wyde, | The rooms and stables spacious were and wide, |
| 29: And wel we weren esed atte beste. | And well we there were eased, and of the best. |
| University of Virginia Library Electronic Text Center. : http://etext.virginia.edu/mideng.browse.html, Cha2Can Available from: Oxford Text Archive no. OTA U-1678-C 1993 [13] | Linkname: Canterbury Tales (prologue) URL: gopher://gopher.vt.edu:10010/02/63/38 |

INSTITUT DE MATHÉMATIQUE PURE ET APPLIQUÉE, UNIVERSITÉ CATHOLIQUE DE LOUVAIN, CHEMIN DU CYCLOTRON,2, B-1348 LOUVAIN-LA-NEUVE (BELGIUM)

*E-mail address*: magnus@anma.ucl.ac.be   , www:   http://math.ucl.ac.be/~magnus/